\documentclass{amsart}
\usepackage[left=2cm, right=2cm, top=2cm]{geometry}
\usepackage[unicode]{hyperref}

\usepackage{latexsym}
\usepackage{indentfirst}
\usepackage{amsxtra}
\usepackage{amssymb}
\usepackage{amsthm}
\usepackage{amsmath}
\usepackage{xcolor}

\usepackage{enumitem}

\newlist{thmlist}{enumerate}{1}
\setlist[thmlist]{label=\arabic*., ref=\thetheor.\arabic*, font=\upshape,noitemsep}

\usepackage{thmtools}

\usepackage{nameref,hyperref}
\usepackage[capitalize]{cleveref}

\newtheorem{theor}{Theorem}
\newtheorem{prop}[theor]{Proposition}
\newtheorem{lemma}[theor]{Lemma}

\theoremstyle{definition}               
\newtheorem{defin}[theor]{Definition}
\newtheorem{ex}[theor]{Example}

\Crefname{thm}{Theorem}{Theorems}
\Crefname{lem}{Lemma}{Lemmas}

\addtotheorempostheadhook[theor]{\crefalias{thmlisti}{theor}}
\addtotheorempostheadhook[lemma]{\crefalias{thmlisti}{lemma}}

\DeclareMathOperator{\End}{End}

\DeclareMathOperator{\id}{id}

\DeclareMathOperator{\Map}{Map}

\setlist[enumerate,1]{label = \arabic*.,ref   = \arabic*, font=\upshape}

\newcommand{\lambdaa}[2]{\lambda_{#1}{#2}}
\newcommand{\rhoo}[2]{\rho_{#1}{#2}}

\title{The algebraic structure of left semi-trusses}
\author{I. Colazzo and A. Van Antwerpen}
\date{}
\address[I. Colazzo]{Department of Mathematics, Vrije Universiteit Brussel, Pleinlaan 2, 1050 Brussel}
\email{ilaria.colazzo@vub.be}
\address[A. Van Antwerpen]{Department of Mathematics, Vrije Universiteit Brussel, Pleinlaan 2, 1050 Brussel}
\email{arne.van.antwerpen@vub.be}
\begin{document}

\maketitle

\begin{abstract}
    The distributive laws of ring theory are fundamental equalities in algebra. However, recently in the study of the Yang-Baxter equation, many algebraic structures with alternative ``distributive" laws were defined. In an effort to study these ``left distributive" laws and the interaction they entail on the algebraic structures, Brzezi\'nski introduced skew left trusses and left semi-trusses. In particular the class of left semi-trusses is very wide, since it contains all rings, associative algebras and distributive lattices. In this paper, we investigate the subclass of left semi-trusses that behave like the algebraic structures that came up in the study of the Yang-Baxter equation. We study the interaction of the operations and what this interaction entails on their respective semigroups. In particular, we prove that in the finite case the additive structure is a completely regular semigroup. Secondly, we apply our results on a particular instance of a left semi-truss called an almost left semi-brace, introduced by Miccoli to study its algebraic structure. In particular, we show that one can associate a left semi-brace to any almost left semi-brace. Furthermore, we show that the set-theoretic solutions of the Yang-Baxter equation originating from almost left semi-braces arise from this correspondence.
\end{abstract}
\section*{Introduction}
The classical examples of a set carrying two operations that interact with each other are rings, associative algebras and lattices. In recent years, several new classes have been defined and actively studied. In the process of studying the Yang-Baxter equation, Rump defined braces in \cite{Ru07}. A left brace, following the equivalent definition of Ced\'o, Jespers and Okni\'nski \cite{CeJO14}, is a set $B$ with an abelian group structure $(B,+)$ and another group structure $(B,\circ)$ such that, for any $a,b,c \in B$, it holds that $a \circ (b+c) = (a \circ b) - a + ( a \circ c)$, where $-a$ denotes the inverse of $a$ in $(B,+)$. It was shown that left braces provide involutive, non-degenerate set-theoretic solutions of the Yang-Baxter equation. Recall, a set-theoretic solution of the Yang-Baxter equation is a tuple $(X,r)$, where $X$ is a set and $r:X \times X \longrightarrow X \times X$ a map such that $$ (r \times \id_X) (\id_X \times r) (r \times \id_X) = (\id_X \times r) (r \times \id_X) (\id_X \times r).$$ Denote $r(x,y) = (\lambda_x(y),\rho_y(x))$. Then, $r$ is called left (resp. right) non-degenerate if $\sigma_x$ (resp. $\gamma_x)$ is bijective for every $x \in X$. The introduction of left braces \cite{Ru07} provided a novel algebraic way to describe all involutive, non-degenerate set-theoretic solutions, which was a problem posed by Drinfel'd \cite{Drinfeld} in 1992. In the wake of this discovery, the notion of a skew left brace was introduced by Guarnieri and Vendramin \cite{GV17}. A skew left brace is a set $B$ equipped with two group operations $(B,+)$ and $(B,\circ)$, where, for any $a,b,c \in B$, it holds that $a \circ (b+c) = (a \circ b) - a + (a \circ c)$. Skew left braces provide bijective non-degenerate set-theoretic solutions of the Yang-Baxter equation and vice versa. Hence, the structure of (skew) left braces has been intensively studied (see, among others, \cite{BaCeJeOk19, CCoSt19, CeGaIvSm18,JeKuVAV19, LeVe16}). In particular, several authors have used ring and group theoretical methods to study these structures \cite{BaCeJeOk19,CSV,JeKuVAV19,JeKuVAV20,konovalov2018skew}.

However, the introduction of these techniques raises several intriguing questions. Can these methods be applied in a completely general algebraic framework? What is the effect of a ``left distributive" law? With these questions in mind Brzezi\'nski introduced skew left trusses \cite{Br19} and left semi-trusses \cite{Br18}. A left semi-truss is a set $B$ equipped with two semigroups $(B,+)$ and $(B,\circ)$ and for every $a \in B$ there exists a map $\lambda_a$ such that for any $b,c \in A$, it holds that $$ a\circ (b+c) = (a \circ b) + \lambda_a(c).$$ If both $(B,+)$ and $(B,\circ)$ are groups, then $(B,+,\circ,\lambda)$ is called a skew left truss. This general framework encompasses the classes of rings, algebras and even distributive lattices. As both skew left trusses and left semi-trusses are concrete objects, an algebraic study is possible and warranted. Recently, Brzezi\'nski \cite{Br20} and  Brzezi\'nski and Rybo{\l}owicz \cite{Br19xa,Br19xb} have introduced and applied a module theory for left trusses in a succesful effort to study both ring and brace theoretical modules in the same framework. 

As most theory and knowledge originates from understanding examples, Brzezi\'nski \cite{Br18} has studied several subclasses of left semi-trusses. In particular, those where the additive semigroup is either left cancellative or an inverse semigroup. In this paper, we continue this study and restrict ourselves to the subclass of brace-like left semi-trusses. In particular, these comprise all (skew) left braces, left semi-braces \cite{CCoSt17,JVA19} and almost left semi-braces \cite{Mi18}. As an application of our results, we study the class of almost left semi-braces, which was introduced by Miccoli \cite{Mi18}. In particular, we show that the additive semi-group of an almost left semi-brace is completely simple and associate to every almost left semi-brace a left semi-brace. Furthermore, as the Yang-Baxter equation is being actively researched, we show that the solution one can associate with an almost left semi-brace is already the associated solution of a left semi-brace. In particular, this shows that this path of generalization will not yield a universal algebraic structure that produces set-theoretic solutions, which is interesting and useful knowledge in its own right. However, several other possibilities for such a structure remain, as evidenced by the recent study of left non-degenerate solutions and their structure monoid by Ced\'o, Jespers and Verwimp \cite{CeJeVe20} and generalized semi-braces by Catino, Colazzo and Stefanelli \cite{CCoSt20}.

\section{Preliminaries}
Brzezi\'{n}ski in \cite{Br19} introduced skew left trusses to deepen the understanding of the nature of the brace-distributive law. 

\begin{defin}
A set $B$ equipped with a group structure $(B,+)$ and group structure $(B,\circ)$ is called a \emph{skew left truss} $(B,+,\circ,\sigma)$, if there exists a map $\sigma:B \longrightarrow B$ such that, for any $a,b,c \in B$, it holds that $$ a \circ (b+c) = \left( a \circ b \right) - \sigma(a) + \left(a \circ b \right),$$ where $-x$ denotes the inverse of $x$ in $(B,+)$.
\end{defin}

\begin{ex}
Every skew left brace $(B,+,\circ)$ is a skew left truss $(B,+,\circ,\id_B)$.
\end{ex}

\begin{ex}
    Let $(G,+)$ be a centerless group. For any $a,b \in G$, denote $a\circ b = a + b - a$. Then, $(G,\circ) \cong (G,+)$ is a group. Furthermore, for any $a,b, c \in G$, it holds that $$ a \circ (b+c) = a + b + c - a = a + b - a + a +c - a = a \circ b + a \circ c.$$ Hence, $(G,+,\circ,\id)$ is a skew left truss.
\end{ex}

In particular, the class of skew left semi-trusses contains all skew left braces. Note that Brzezi\'nski has extended his notion of skew left trusses such that the additive group $(B,+)$ is replaced by a heap \cite{Br19}. For the continuity of the exposition, we chose to present the definition above.

Later, Brzezi\'nski \cite{Br18} extended this to the notion of left semi-trusses. 

\begin{defin}
Let $B$ be a set with two associative operations $(B,+)$ and $(B,\circ)$. If, for any $a \in B$, there exists a map $\lambda_a: B \longrightarrow B$ such that for any $b,c \in B$, it holds that $$ a \circ (b+ c) = \left( a \circ b\right) + \lambda_a(b),$$
then $(B,+,\circ,\lambda)$ is called a left semi-truss.
\end{defin}

Clearly, the class of semi-trusses contains all skew left trusses, rings, associative algebras and distributive lattices. This entails that it will prove difficult to present deep results on this class. However, one may examine large subclasses in an effort to gain a deeper understanding of the larger phenomenon. In particular, Brzezi\'nski \cite{Br18} focused on semi-trusses with $\left(B,+\right)$ a left cancellative or inverse semigroup. He showed that in these cases, for many semigroups, the ``distributive" law can be rewritten in a fashion that closely resembles the ``distributive" law of skew left braces.

In this paper, we examine the effect of the ``distributive" law on the interaction between the additive and multiplicative structures of the left semi-truss. The following example shows that one will have to restrict to interesting subclasses. In particular, it shows that every semigroup can occur as the additive or multiplicative semigroup of a left semi-truss.
\begin{ex}
    Let $(S,+)$ be any semigroup. Then, $(S,+,+,\id_S)$ is a left semi-truss. 
\end{ex}

\begin{ex}
    Let $(B,+)$ be a group. For any $a,b \in B$, denote $a\circ b = a+b-a$. Then, $(B,+,\circ,\id)$ is a left semi-truss.
\end{ex}

We examine more closely the class of left semi-trusses that resembles (skew) left braces.

\begin{defin}\label{def:semitrusses}
	Let $B$ be a set with two operations $+$ and $\circ$ such that $\left(B,+\right)$ is a semigroup, $\left(B,\circ\right)$ a group, and  $\lambda:B\to \End\left(B,+\right)$ a morphism from $\left(B,\circ\right)$ into the endomorph of $\left(B,+\right)$.
	We say that $\left(B,+,\circ, \lambda\right)$ is a \emph{brace-like left semi-truss} if, for all $a,b, c \in B$, the following condition holds
	\begin{align}\label{eq:semitrusses}
	a \circ\left(b+c\right)=a \circ b + \lambdaa{a}{\left(c\right)}.
	\end{align}
	We call $\left(B,+\right)$ the \emph{additive semigroup}, and $\left(B,\circ\right)$ the \emph{multiplicative group} of the brace-like left semi-truss $\left(B,+,\circ,\lambda\right)$. 
\end{defin}

The class we restrict to contains all left semi-braces \cite{CCoSt17,JVA19}. 
Recall that a \emph{(left cancellative) left semi-brace} is a triple $\left(B,+,\circ\right)$ such that $\left(B,+\right)$ is a (left cancellative) semigroup, $\left(B,\circ\right)$ is a group and, for all $a,b,c \in B$, the following condition holds
\begin{align*}
a \circ\left(b+c\right) = a \circ b + a \circ\left(\overline{a}+ c\right),
\end{align*}
where $\overline{a}$ denotes the inverse of $a$ in $\left(B,\circ\right)$.
If $(B, +,\circ)$ is a left semi-brace, then it is a brace-like left semi-truss with $\lambdaa{a}{\left(b\right)} = a\circ\left(\bar{a}+b\right)$, for every $a,b \in  B$, since by \cite[Proposition 3]{CCoSt17} and \cite[Lemma 2.12]{JVA19} the function $\lambda: B\to \End{\left(B,+\right)}, a \to \lambdaa{a}{}$ is a morphism. Vice versa if $\left(B,+,\circ,\lambda\right)$ is a brace-like left semi-truss with $\lambdaa{a}{\left(b\right)} = a\circ\left(\bar{a}+b\right)$, then $\left(B,+,\circ\right)$ is a left semi-brace.	

Let $\left(B,+,\circ,\lambda\right)$ be a brace-like left semi-truss. 
The element $1$ denotes the identity of the multiplicative group $\left(B,\circ\right)$.

First, we should remark that one will not be able to give restrictions on the properties of the group $(B,\circ)$ of a brace-like semi-truss.

\begin{ex}
    Let $(B,\circ)$ be any group. Denote for any $a,b \in B$, the operation $a+b = a$, i.e. $(B,+)$ is a left zero semigroup. Then, consider any semigroupmorphism $\lambda: (B,\circ) \longrightarrow \Map(B,B)$. Then, $(B,+,\circ,\lambda)$ is a brace-like left semi-truss.
\end{ex}

However, in the next Section we shall prove that there does exist a restriction on the additive semigroup $(B,+)$.

\section{The algebraic structure of a semi-truss}\label{sec:algstructure}

In this section, we focus on the additive structure of a brace-like left semi-truss $B$. In particular, we prove that if $B$ is finite, then the additive semigroup $\left(B,+\right)$ is a completely simple semigroup.

\begin{lemma}\label{lemma:one}
	Let $\left(B,+,\circ, \lambda\right)$ be a brace-like left semi-truss and $1$ be the identity of $\left(B,\circ\right)$. The following statements hold:
	\begin{thmlist}
		\item $\lambdaa{1}{}\lambdaa{1}{}=\lambdaa{1}{}$\label{lemma:one1};
		\item $\lambdaa{1}{\lambdaa{a}{}}=\lambdaa{a}{}=\lambdaa{a}{\lambdaa{1}{}}$, for any $a \in B$;\label{lemma:one2}
		\item $a+b=a +\lambdaa{1}{\left(b\right)}$, for all $a,b \in B$. In particular $B+B=B+\lambdaa{1}{\left(B\right)}$; \label{lemma:one3}
		\item $1+B$ is a subsemigroup of $\left(B,\circ\right)$.
	\end{thmlist}
\begin{proof}
	Since $\lambda$ is a morphism and $\left(B,\circ\right)$ is a group with identity $1$, both $1.$ and $2.$ easily follow. 
		\begin{enumerate}[label={\arabic*}.]
		\item[3.] Let $a,b\in B$. Then 
		\begin{align*}
		a+b= 1 \circ\left(a+b\right) = 1\circ a + \lambdaa{1}{\left(b\right)} = a +\lambdaa{1}{\left(b\right)}.
		\end{align*}
		\item[4.] Let $a,b \in B$. Then
		\begin{align*}
		\left(1+a\right)\circ \left(1+b\right) = \left(1+a\right)\circ 1 + \lambdaa{\left(1+a\right)}{\left(b\right)} = 1 + a + \lambdaa{\left(1+a\right)}{\left(b\right)}\in 1+B.
		\end{align*}
		Therefore, $1+B$ is a subsemigroup of $\left(B,\circ\right)$. 
		\end{enumerate}
\end{proof}
\end{lemma}
\begin{lemma}\label{lemma:zeroelement}
 	Let $\left(B,+,\circ, \lambda\right)$ be a brace-like left semi-truss. If $B$ has at least two elements then $\left(B,+\right)$ does not contain a zero element.
 	\begin{proof}
		Suppose $\theta$ is a zero element of $\left(B,+\right)$. Then, for all $a,b \in B$, it holds that
	\begin{align*}
a \circ \theta = a \circ \left(\bar{a}\circ \theta + \theta\right) = a \circ \bar{a}\circ \theta + \lambdaa{a}{\left(\theta\right)} = \theta+\lambdaa{a}{\left(\theta\right)}=\theta.
\end{align*}
As $\left(B,\circ\right)$ is a group, it follows that $B=\left\{\theta\right\}$, a contradiction.
 	\end{proof}
\end{lemma}
\begin{lemma}\label{lemma:B=B+B}
	Let $\left(B,+,\circ, \lambda\right)$ be a brace-like left semi-truss and $1$ be the identity of $\left(B,\circ\right)$. The following statements hold:
	\begin{thmlist}
		\item $1 \in B+B$ and $B+B$ is a subgroup of the multiplicative group $\left(B,\circ\right)$;\label{lemma:B=B+B1}
		\item $\left(B\setminus\left(B+B\right)\right)\cup \left\{1\right\}$ is a subgroup of $\left(B,\circ\right)$;\label{lemma:B=B+B2}
		\item $B=B+B$.\label{lemma:B=B+B3}
	\end{thmlist}
\begin{proof}\mbox{}
	\begin{enumerate}[label={\arabic*}.]
		\item Let $a,b \in B$. Then
		\begin{align*}
			1=\overline{a+b}\circ\left(a+b\right) = \overline{a+b}\circ a + \lambdaa{\overline{a+b}}{\left(b\right)}\in B+B.
		\end{align*}
		Moreover, let $b,c \in B$. Then
		\begin{align*}
			\overline{a+b}\circ\left(c+d\right) = \overline{a+b}\circ c + \lambdaa{\overline{a+b}}{\left(d\right)}\in B+B.
		\end{align*}
		Hence $\left(B+B,\circ\right)$ is a subgroup of $\left(B,\circ\right)$.
		\item Put $D:=\left(B\setminus \left(B+B\right)\right)\cup \left\{1\right\}$. We need to show that $\overline{a}\circ b\in D$, for all $a,b\in D$. Suppose the contrary, i.e., assume there exist distinct $a,b \in D$ such that $\overline{a}\circ b=s+t$, for some $s,t\in B$. Hence
		\begin{align*}
			b = a \circ\left(s+t\right) = a\circ s + \lambdaa{a}{\left(t\right)}\in B+B,
		\end{align*}
		which is a contradiction.
		\item The claim is trivially true if $B=\left\{1\right\}$. Suppose $\left|B\right|\geq 2$. From $1.$ and $2.$, we know that the group $\left(B,\circ\right)$ is the union of two subgroups $B+B$ and  $\left(B\setminus \left(B+B\right)\right)\cup \left\{1\right\}$. Because the intersection of these groups is $\left\{1\right\}$ it follows that $B= B+B$ or $B+B = \left\{1\right\}$. In the second case, $b+1 = 1 = 1+b$, for any $b\in B$, in contradiction with \cref{lemma:zeroelement}, since $\left| B\right| \geq 2$.
	\end{enumerate}
\end{proof}
\end{lemma}

\begin{lemma}\label{lemma:idempotentelement}
	Let $\left(B,+,\circ, \lambda\right)$ be a brace-like left semi-truss with $B$ finite. Then there exists an idempotent element $z$ of $\left(B,+\right)$ such that $z+B+z = 1+B+z$ is a subgroup of $\left(B,+\right)$.
	\begin{proof}
		For any $a \in B$, define $a^1 := a$, for every integer $n\geq 1$, and $a^{n+1} =a \circ a^n$.
		
		First, suppose that for all $b\in B$, $1+b=1$. Then, in particular, $1+1=1$, i.e., $1$ is an idempotent element. Moreover, $1+B+1=\left\{1\right\}$ and the claim is trivially true. 
		
		Otherwise, since $B$ is finite, for any $b \in B$ such that $1+b\neq 1$, there exists an integer $n_b\geq 2$ such that $\left(1+b\right)^{n_b} = 1$. Then
		\begin{align*}
			1 &= \left(1+b\right)^{n_b-1} \circ \left(1+b\right)\\
			&= \left(1+b\right)^{n_b-1} + \lambdaa{\left(1+b\right)^{n_b-1}}{\left(b\right)}.\end{align*}
		Continuing, we obtain inductively that 
			$$1=1+b + \sum_{i=1}^{n_b-1}\lambdaa{\left(1+b\right)^{i}}{\left(b\right)}.$$
		Set $z_b:=b + \sum_{i=1}^{n_b-1}\lambdaa{\left(1+b\right)^{i}}{\left(b\right)}$. Then $1+z_b=1$. 
		
		First, if $z_b=1$, then $1$ is idempotent and it is the identity of $1+B+1$. Moreover, since $B$ is finite, for any $a\in B$ such that $1+a+1\neq 1$, there exists a integer $k\geq 2$, such that $\left(1+a+1\right)^k=1$. Set $\underline{a}:=\sum_{i=1}^{k-1}\lambdaa{\left(1+a+1\right)^i}{\left(a+1\right)}$, then
		\begin{align*}
		    \left(1+a+1\right)+\left(1+\underline{a}+1\right) 
		    &= 1+a+1+\sum_{i=1}^{k-1}\lambdaa{\left(1+a+1\right)^i}{\left(a+1\right)}+1\\
		    &=\left(1+a+1\right)^{k}+1=1+1=1,
		\end{align*}
		i.e., every element in $1+B+1$ has a right inverse. Therefore $1+B+1$ is a group. 
		Otherwise, $z_b \neq 1$. Since $B$ is finite, there exists an integer $m_b\geq 2$, such that $z_b^{m_b}=1$ and, so, 
		\begin{align*}
			\lambda_{z_b}^{m_b}=\underbrace{\lambdaa{z_b}{\lambdaa{z_b}{\ldots \lambdaa{z_b}{}}}}_{m_b\mbox{ times}}=\lambdaa{1}{}.
		\end{align*}	
		Moreover, note that $z_b = z_b\circ 1 = z_b\circ \left(1+z_b\right)= z_b + \lambdaa{z_b}{\left(z_b\right)}$, then
		$$z_b = z_b+\lambdaa{z_b}{\left(z_b\right)} = z_b +\lambdaa{z_b}{\left(z_b+\lambdaa{z_b}{z_b}\right)} = z_b + \lambdaa{z_b}{\left(z_b\right)}+\lambda^2_{z_b}\left(z_b\right).$$
			If we continue this inductively, we obtain that 
		    $$z_b ={z_b} + \lambda_{z_b}\left({z_b}\right)+\cdots + \lambda_{z_b}^{m_b-1}\left(z\right)+\lambda^{m_b}_{z_b}\left({z_b}\right)+ \cdots+\lambda_{z_b}^{2m_b-1}\left(z_b\right).$$
		Set $\underline{z_b}:=z_b + \lambda_{z_b}\left(z_b\right)+\cdots + \lambda_{z_b}^{m_b-1}\left({z_b}\right)$, then
		\begin{align*}
			z_b=\underline{z_b}=\underline{z_b}+\lambda_{z_b}^{m_b}\left(\underline{z_b}\right)=\underline{z_b}+\lambdaa{1}{\left(\underline{z_b}\right)} =\underline{z_b}+\underline{z_b}=z_b+z_b.
		\end{align*}
		In particular, this shows that $z_b$ is an idempotent. From now on, we focus on $b \in 1+B$. Then, by definition, there exists $z'_b$ such that $z_b=1+z'_b$.
		Then
		\begin{align*}
			1 &= z_b^{m_b} = \left(1+z'_b\right)^{m_b} = 1+z'_b+\sum_{i=1}^{m_b-1}\lambdaa{z_b^i}{\left(z'_b\right)} 
			= z_b +\sum_{i=1}^{m_b-1}\lambdaa{z_b^i}{\left(z'_b\right)} \\
			&=z_b+\left(z_b +\sum_{i=1}^{m_b-1}\lambdaa{z_b^i}{\left(z'_b\right)} \right)
			=z_b + z_b^{m_b}\\
			&=z_b +1.
		\end{align*}

		Finally, we prove that $1+B+z_b$, with $b\in 1+B$,  is a subgroup of $\left(B,+\right)$ with identity $z_b$. First, $$z_b = z_b+z_b = 1+ z'_b+ z_b \in 1+B+z_b.$$ 
		Moreover, for all $a\in B$,
		\begin{align*}
			\left(1+ a + z_b\right) + z_b = 1 + a + z_b
		\end{align*}
		and
		\begin{align*}
			z_b+\left(1+ a + z_b\right) = \left(z_b+1\right) + a + z_b = 1 + a + z_b.
		\end{align*}
		Furthermore, since $B$ is finite for all $a\in B$, there exists a positive integer $k$ such that $\left(1+a+z_b+1\right)^k=1$. If $k=1$, then $\left(1+a+z_b\right)+\left(1+z'_b+z_b\right) = 1+z'_b+z_b =z_b+z_b=z_b$, otherwise
		\begin{align*}
			1 = \left(1+a+z_b+1\right)^k = 1+a+z_b+1 + \sum_{i=1}^{k-1}\lambdaa{\left(1+a+z_b+1\right)^i}{\left(a+z_b+1\right)},
		\end{align*}
		set $\underline{a}:= \sum_{i=1}^{k-1}\lambdaa{\left(1+a+z_b+1\right)^i}{\left(a+z_b+1\right)}$, then
		\begin{align*}
			\left(1+a+z_b\right)+\left(1+\underline{a}+z'_b+z_b\right) =  \left(1+a+z_b+1+\underline{a}\right)+z'_b+z_b = 1+z'_b+z_b=z_b+z_b=z_b, 
		\end{align*}
		i.e., every element in $1+B+z_b$ has a right inverse. Therefore $1+B+z_b$ is a group. 
		Finally, for any $a\in B$, $1+a+z_b= z_b+1+a+z_b \in z_b+B+z_b$ and $z_b+a+z_b= 1+z'_b+a+z_b\in 1+B+z$. Hence, $z_b+B+z_b=1+B+z_b$.
	\end{proof}
\end{lemma}

\begin{lemma}\label{lemma:B+z+B}
		Let $\left(B,+,\circ, \lambda\right)$ be a brace-like left semi-truss. If there exists an idempotent element $z$ of $\left(B,+\right)$ such that $z+B+z$ is a subgroup of $\left(B,+\right)$. Then, for any $c\in B$, $c\in c+z+B$ and, in particular, $B+z+B=B$.
		\begin{proof}
			Let $c \in B$.
			If $d\in B$ is such that $z+d+z$ is the inverse of $z+\lambdaa{z}{\lambdaa{\bar{c}}{\left(z\right)}}+z$ in $z+B+z$ and $b := \lambdaa{c}{\lambdaa{\bar{z}}{\left(z+d+z\right)}}$. Then
			$c=c+z+b$. Indeed
			\begin{align*}
				c &= c\circ \bar{z}\circ z 
				=c \circ \bar{z} \circ \left(z+\lambdaa{z}{\lambdaa{\bar{c}}{\left(z\right)}}+z+z+d+z \right)\\
				&=c\circ\bar{z}\circ\left(z+\lambdaa{z}{\lambdaa{\bar{c}}{\left(z\right)}}+\lambdaa{z}{\lambdaa{\bar{c}}{\left(\lambdaa{c}{\lambdaa{\bar{z}}{\left(z+d+z\right)}}\right)}}\right)\\
				&=c\circ\bar{z}\circ\left(z + \lambdaa{z}{\lambdaa{\bar{c}}{\left(z+b\right)}}\right)\\
				&=c\circ\bar{z}\circ z \circ \bar{c}\circ\left(c+z+b\right)\\
				&=c+z+b\in B+z+B.
			\end{align*}
			Therefore $B=B+z+B$.
		\end{proof}
\end{lemma}

\begin{theor}\label{ther:completelysimple}
		Let $\left(B,+,\circ, \lambda\right)$ be a brace-like left semi-truss. If there exists an idempotent element $z$ of $\left(B,+\right)$ such that $z+B+z$ is a subgroup of $\left(B,+\right)$. Then $z$ is a primitive idempotent and $\left(B,+\right)$ is completely simple semigroup. In particular, this holds if $B$ is finite.
		\begin{proof}
			Let $b\in B$. Since $z+B+z$ is a group, there exists an element $\underline{b}\in B$ such that $z+b+z+z+\underline{b}+z = z$. Moreover, by \cref{lemma:B+z+B}, for any $c\in B$, there exists an element $d\in B$ such that $c=c+z+d$. Then  
			\begin{align*}
				c = c +z+ d = \left(c +z\right)+b+\left(z+z+\underline{b}+z+ d\right) \in B+b+B.
			\end{align*}
			Hence, $B=B+b+B$. Therefore, every principal ideal, and thus every ideal of $B$ is trivial, i.e., $\left(B,+\right)$ is a simple semigroup. Since $z+B+z$ is a subgroup of $\left(B,+\right)$, then $z$ is a primitive idempotent of $\left(B,+\right)$. Hence $\left(B,+\right)$ is a completely simple semigroup.
		\end{proof}
\end{theor}

 Theorem \ref{ther:completelysimple} shows that the role of idempotents of $(B,+)$ should not be underestimated. Hence, it is an interesting question whether these idempotents form a subsemi-truss, like in the case for left semi-braces \cite{CCoSt17,JVA19}. Denote the set of idempotents of $(B,+)$ by $E\left(B\right)$. The following lemma proves that $E\left(B\right)$ is closed under any $\lambda_a$, for each left semi-truss $\left(B,+,\circ,\lambda\right)$.

\begin{lemma}\label{lemma:lambda}
	Let $\left(B,+,\circ,\lambda\right)$ be a left semi-truss. For any $a \in B$, $\lambda_a\left(E\left(B\right)\right)\subseteq E\left(B\right)$.
	\begin{proof}
		If $e \in E\left(B\right)$, then 
		$
		\lambda_a\left(e\right) = \lambda_a\left(e+e\right) = \lambda_a\left(e\right)+\lambda_a\left(e\right).
		$
		Hence, $\lambda_a\left(e\right) \in E\left(B\right)$. 
	\end{proof}
\end{lemma}

However, in general $E\left(B\right)$ is not closed with respect to $\circ$ as shown by the following example.

\begin{ex}
    Let $B$ be the Klein group defined by
	\begin{align*}
	B:=\left\langle\left. a,b\ \right|\ a^2=b^2=\left(ab\right)^2 = 1\right\rangle.
	\end{align*}
	Every element $x \in B$ can be represented by $x = a^ib^j$ with $i,j \in \mathbb{Z}$. Clearly, two elements $a^ib^j$ and $a^hb^k$ are equal if and only if $i \equiv h \mod 2$ and $j \equiv k \mod 2$. Define
	\begin{align*}
	a^ib^j+a^hb^k := a^{1+i+h}b^{j}
	\end{align*}
	and 
	\begin{align*}
	    \lambda: B \longrightarrow \End\left(B,+\right), \quad a^{i}b^{j}\longmapsto f
	\end{align*}
	where $f:B\to B,\ a^{h}b^{k}\mapsto a^{h}$. Then $\left(B,+,\circ,\lambda\right)$ is a left semi-truss. Indeed
	\begin{align*}
	    \lambdaa{a^ib^j}{\lambdaa{a^hb^k}{\left(a^yb^z\right)}} = ff\left(a^yb^z\right) = f\left(a^y\right) = a^y = f\left(a^yb^z\right)=\lambdaa{a^ib^j\circ a^hb^k}{\left(a^yb^z\right)},
	\end{align*}
	\begin{align*}
	    \lambdaa{a^ib^j}{\left(a^hb^k+a^yb^z\right)} = \lambdaa{a^ia^j}{\left(a^{1+h+y}b^k\right)} =a^{1+h+y}
	\end{align*}
	and 
	\begin{align*}
	    \lambdaa{a^ib^j}{\left(a^hb^k\right)} + \lambdaa{a^ib^j}{\left(a^yb^z\right)}
	    =a^h+a^y = a^{1+h+y}.
	\end{align*}
	Moreover
	\begin{align*}
	    a^ib^j\circ\left(a^hb^k+a^yb^z\right) = a^ib^j\circ\left(a^{1+h+y}b^k\right) 
	    = a^{i+1+h+y}b^{j+k}
	\end{align*}
	and
	\begin{align*}
	    a^ib^j\circ a^hb^k + \lambdaa{a^ib^j}{\left(a^yb^z\right)}
	    =a^{i+h}b^{j+k} + a^y = a^{1+i+h+y}b^{j+k}.
	\end{align*}
	Finally, $a+a = a$, i.e., $a\in E\left(B\right)$, but $a\circ a =1$ and $1+1=a$, i.e., $a\circ a \notin E\left(B\right).$
\end{ex}

In particular, this shows that in general the additive and multiplicative semigroups of a brace-like left semi-truss fit together very differently than those of a left semi-brace, where a nice decomposition theorem \cite{CCoSt17,JVA19} can be obtained.

\section{Almost left semi-brace and solutions}\label{sec:amostsemi}

In this section, we focus on a particular instance of brace-like left semi-trusses called almost left semi-braces. Our definition extends the one introduced by Miccoli \cite{Mi18}, where one works under the restriction that the semigroup $\left(B, +\right)$ is left cancellative and, for all $a,b \in B$, it holds
\begin{equation}\label{eq:defAlmost2}
	\left(\iota\left(a\right)+b\right)\circ \iota\left(1\right) = \iota\left(a\right)+b\circ\iota\left(1\right),
\end{equation}
where $1$ is the identity of the group $\left(B,\circ\right)$.

\begin{defin}\label{def:almost_semibrace}
	Let $B$ be a set with two operations $+$ and $\circ$ such that $\left(B,+\right)$ is a semigroup, $\left(B,\circ\right)$ is a group, and $\iota:B \to B$ is a map such that, for all $a,b \in B$, the following condition holds
	\begin{equation}\label{eq:defAlmost1}
	    \iota\left(a\circ b\right) = \bar{b}\circ \iota\left(a\right),
	\end{equation}
	where $\bar{b}$ denotes the inverse of $b$ in $\left(B, \circ\right)$.
	We say that $\left(B, +, \circ, \iota\right)$ is an \emph{almost left semi-brace} if, for all $a,b,c\in B$, it holds
	\begin{align}
	\label{eq:defAlmost3}a \circ \left(b + c\right) = a \circ b + a \circ \left(\iota\left(a\right) + c\right).
	\end{align}
	We call $\left(B,+\right)$ the \emph{additive semigroup}, and $\left(B,\circ\right)$ the \emph{multiplicative group} of the almost left semi-brace $\left(B,+,\circ,\iota\right)$.
\end{defin}

If $\left(B,+,\circ\right)$ is a left semi-brace, then it is an almost left semi-brace with $\iota:B\to B$ defined by $\iota\left(a\right)=\bar{a}$, for any $a\in B$. Conversely, if $\left(B,+,\circ,\iota\right)$ is an almost semi-brace such that $\iota\left(a\right)= \bar{a}$, then $\left(B,+,\circ\right)$ is a left semi-brace.

The following lemma proves that given an almost left semi-brace $\left(B,+,\circ,\iota\right)$, then it is possible to define a morphism $\lambda$ such that $\left(B,+\circ,\lambda\right)$ is a brace-like left semi-truss.

\begin{lemma}\label{lemma:lambda_almost}
	Let $\left(B,+,\circ,\iota\right)$ be an almost left semi-brace. For any $a \in B$, it follows that $\lambda_a \in \End\left(B,+\right)$, where $\lambda_a\left(b\right):=a \circ \left(\iota\left(a\right)+b\right)$, for every $b \in B$.
	Moreover, $\lambda: \left(B,\circ\right) \to \End\left(B,+\right), a\mapsto \lambda_a$ is a homomorphism and, for any $a\in B$, it holds that $\lambda_a\left(E\left(B\right)\right)\subseteq E\left(\iota\left(1\right)+B\right)$.
	\begin{proof}
		First, note that defining, for all $a,b\in B$, $\lambdaa{a}{\left(b\right)}:=a\circ\left(\iota\left(a\right)+b\right)$, \eqref{eq:defAlmost3} is clearly equivalent to
		\begin{align}\label{eq:defAlmost3'}
		a\circ \left(b+c\right) = a \circ b + \lambda_a\left(c\right).
		\end{align}
		Let $a, b, x, y\in B$. By \eqref{eq:defAlmost3'} it holds that
		\begin{align*}
		\lambda_a\left(x+y\right) &= a \circ \left(\iota\left(a\right)+x+y \right)\\
		&= a\circ \left(\iota\left(a\right)+x\right) + \lambda_a\left(y\right)\\
		&=\lambda_a\left(x\right)+\lambda_a\left(y\right),
		\end{align*}
		and, by \eqref{eq:defAlmost1} and \eqref{eq:defAlmost3'}, it follows that
		\begin{align*}
		\lambda_{a\circ b}\left(x\right) &= a \circ b\circ \left(\iota\left(a\circ b\right)+x\right)\\
		&=a\circ b\circ \left(\overline{b}\circ \iota\left(a\right)+x\right)\\
		&=a \circ \left(b \circ \bar{b}\circ \iota\left(a\right) + \lambda_b\left(x\right)\right)\\
		&=a \circ \left( \iota\left(a\right) +\lambda_b\left(x\right)\right)\\
		&=\lambda_a\lambda_b\left(x\right).
		\end{align*}
		Moreover, if $e \in E\left(B\right)$, then 
		$
		\lambda_a\left(e\right) = \lambda_a\left(e+e\right) = \lambda_a\left(e\right)+\lambda_a\left(e\right).
		$
		Hence, $\lambda_a\left(e\right) \in E\left(B\right)$. Furthermore,
		\begin{align*}
		\iota\left(1\right)+\lambda_a\left(e\right) = a \circ \bar{a}\circ \iota\left(1\right) + \lambda_a\left(e\right)
		=a\circ \iota\left(a\right) +\lambda_a\left(e\right)
		=a \circ\left(\iota\left(a\right)+e \right) = \lambda_a\left(e\right).
		\end{align*}
		Therefore, $\lambda_a\left(E\left(B\right)\right) \subseteq E\left(\iota\left(1\right)+B\right)$.
	\end{proof}
\end{lemma}

Hence, if $\left(B,+,\circ,\iota\right)$ is an almost left semi-brace, defining $\lambda:\left(B,\circ\right)\to \End\left(B,+\right)$ as in \cref{lemma:lambda}, then $\left(B,+,\circ,\lambda\right)$ is a brace-like left semi-truss. But not all brace-like left semi-trusses are almost left semi-braces. Just consider the following example to be convinced.

\begin{ex}
	Let $\left(B,\circ\right)$ be a group, $f:B\to B$ a map such that $f^2=f$ and denote $\left(B,+\right)$ the left zero semigroup on the set $B$. Then $\left(B,+,\circ, \lambda\right)$ where $\lambda:B\to \End\left(B,+\right), a \mapsto f$ is a brace-like left semi-truss, but, in general, is not an almost semi-brace. For instance, suppose that $f=\id_B$ and there exists $\iota:B\to B$ such $\lambdaa{a}{\left(b\right)}=a\circ\left(\iota\left(a\right)+b\right)$. Hence, $b=\lambdaa{a}{\left(b\right)} =a\circ\iota\left(a\right)$, for any $a,b\in B$. In particular, for $b=a$, we get $\iota\left(a\right)=1$. Therefore, $b =\lambdaa{a}{\left(b\right)}=a$, for any $a,b\in B$, a contradiction, if $\left|B\right|>1$.
\end{ex}

If $\left(B_1,+_1,\circ_1,\iota_1\right)$ and $\left(B_2,+_2,\circ_2,\iota_2\right)$ are almost left semi-braces then a map $f : B_1 \to B_2$ is a homomorphism of almost left semi-braces if $f$ is a semigroup homomorphism from $\left(B_1, +_1\right)$ to $\left(B_2, +_2\right)$, $f$ is a group homomorphism from $\left(B_1,\circ_1\right)$ to $\left(B_2,\circ_2\right) $ and, $f\iota_1 = \iota_2f$.
Hence, a left semi-brace $\left(B,+,\circ\right)$ considered as an almost left semi-brace can not be isomorphic to an almost left semi-brace $\left(B,+,\circ,\iota_B\right)$ with $\iota_B\left(1\right)\neq1$. Indeed, such an isomorphism $f$ has to satisfy $\iota_B\left(1\right) = \iota_Bf\left(1\right) = f(1) = 1$. 
However, we can associate a left semi-brace to any almost left semi-brace. To this purpose, we make the following preliminary result. 

\begin{lemma}\label{lemma:24}
	Let $\left(B,+,\circ,\iota \right)$ be an almost left semi-brace. The following properties hold.
	\begin{thmlist}
		\item $a+b = a+\lambdaa{1}{\left(b\right)} = a+ \iota\left(1\right) + b$, for all $a,b\in B$. In particular $B+B = B + \iota\left(1\right)+B$\label{lemma:24.1}
		\item $\iota$ is bijective and $\iota\left(a\right) = \bar{a}\circ\iota\left(1\right)$, for every $a\in B$.\label{lemma:24.2}
	\end{thmlist}
	\begin{proof}\mbox{ }
		\begin{enumerate}[label={\arabic*}.]
			\item Let $a,b \in B$. Then, by \cref{lemma:one3} and \cref{lemma:lambda_almost}, $a+b=a+\lambdaa{1}{\left(b\right)}$ and 
				\begin{align*}
					a+b = a+\lambdaa{1}{\left(b\right)}= a + 1 \circ\left(\iota\left(1\right)+ b\right) =a + \iota\left(1\right)+b.
				\end{align*}
			\item Let $a \in B$. By \eqref{eq:defAlmost1} in \cref{def:almost_semibrace}, we get
				\begin{align}\label{eq:iotalemma}
					\iota\left(a\right) = \iota\left(1\circ a\right) = \bar{a}\circ\iota\left(1\right).
				\end{align}
				Furthermore, $\iota$ is bijective. Indeed, if $a,b\in B$, such that $\iota\left(a\right)=\iota\left(b\right)$, then by \eqref{eq:iotalemma} $\bar{a}\circ\iota\left(1\right) =\bar{b}\circ\iota\left(1\right)$ and so $a=b$. Furthermore, if $b \in B$, then by \eqref{eq:iotalemma},
				\begin{align*}
						\iota\left(\iota\left(1\right)\circ \bar{b}\right) = b \circ \overline{\iota\left(1\right)}\circ\iota\left(1\right) = b \circ 1 = b.
				\end{align*}
		\end{enumerate}
	\end{proof}
\end{lemma}
Applying Theorem \ref{ther:completelysimple} to almost left semi-braces, we obtain the following.

\begin{prop}
    Let $(B,+,\circ,\iota)$ be a finite almost left semi-brace. Then, $(B,+)$ is a completely simple semigroup, where $1$ is a primitive idempotent.
\end{prop}

Now, let $\left(B,+,\circ,\iota\right)$ be an almost left semi-brace. Define the following operation on $B$:
\begin{align*}
\oplus: B \times B \longrightarrow B, \quad \left(a,b\right) \longmapsto \iota^{-1}\left(\iota\left(a\right)+\iota\left(b\right)\right)
\end{align*}

Note that, from \eqref{eq:defAlmost1}, if $\left(B,+,\circ,\iota\right)$ is an almost left semi-brace, then it follows that 
\begin{align}\label{eq:defAlmost2inv}
\iota^{-1}\left(a\circ b\right) = \iota^{-1}\left(b\right)\circ \bar{a}
\end{align}
for all $a,b \in B$. 

\begin{prop}
	If $\left(B,+,\circ,\iota\right)$ is an almost left semi-brace and $\oplus$ defined as above, then $\left(B,\oplus,\circ^{op}\right)$ is a left semi-brace, called \emph{left semi-brace associated to the almost left semi-brace $\left(B,+,\circ,\iota\right)$}.
	\begin{proof}
		Clearly $\left(B,\circ^{op}\right)$ is a group. Moreover, if $a,b,c, \in B$, then 
		\begin{align*}
		\left(a\oplus b\right) \oplus c = \iota^{-1}\left(\iota\left(a\right) +\iota\left(b\right)\right) \oplus c = \iota^{-1}\left(\iota\left(a\right) + \iota\left(b\right) + \iota\left(c\right)\right)
		\end{align*}
		and
		\begin{align*}
		a \oplus\left(b \oplus c\right) = a \oplus\iota^{-1}\left(\iota\left(b\right) +\iota\left(c\right)\right)
		= \iota^{-1}\left(\iota\left(a\right) + \iota\left(b\right) + \iota\left(c\right)\right).
		\end{align*}
		Thus $\left(B,\oplus\right)$ is a semigroup. Finally, if $a,b,c\in B$ then
		\begin{align*}
		a\circ^{op}\left(b\oplus c\right) = \iota^{-1}\left(\iota\left(b\right) + \iota\left(c\right)\right) \circ a
		\end{align*}
		and by \eqref{eq:defAlmost2inv} and \eqref{eq:defAlmost3}
		\begin{align*}
		a \circ^{op} b \oplus a \circ^{op} \left(\bar{a} \oplus c\right)
		&= b \circ a \oplus \left(\bar{a} \oplus c\right) \circ a
		= b\circ a \oplus \iota^{-1}\left(\iota\left(\bar{a}\right) + \iota\left(c\right)\right)\circ a\\
		&=b \circ a \oplus \iota^{-1}\left(\bar{a}\circ \left(\iota\left(a\right) + \iota\left(c\right)\right)\right)\\
		&= \iota^{-1}\left(\iota\left(b\circ a\right) + \bar{a}\circ \left(\iota\left(a\right) +\iota\left(c\right)\right)\right)\\
		&= \iota^{-1}\left(\bar{a}\circ\iota\left(b\right) + \bar{a}\circ \left(\iota\left(a\right) +\iota\left(c\right)\right)\right)\\
		&=\iota^{-1}\left(\bar{a}\circ\left(\iota\left(b\right) +\iota\left(c\right)\right)\right)\\
		&=\iota^{-1}\left(\left(\iota\left(b\right) +\iota\left(c\right)\right)\right)\circ a.
		\end{align*}
		Hence, $\left(B,\oplus,\circ^{op}\right)$ is a left semi-brace.
	\end{proof}
\end{prop}

As in \cite[Corollary 2.7]{Br18} we can associate a solution with an almost left semi-brace via its associated left semi-brace, whenever it is possible to associate a solution to such a semi-brace. But, under some assumptions, we can associate a solution directly to an almost semi-brace. In the following, we prove that these solutions are isomorphic.

Following the idea in \cite{CCoSt19x}, we provide a necessary and sufficient condition for obtaining a solution from an almost left semi-brace.

\begin{theor}\label{theor:solutionNew}
	Let $\left(B,+,\circ, \iota\right)$ be an almost left semi-brace such that \eqref{eq:defAlmost2} is satisfied, for all $a,b\in B$. The map $r_B:B\times B\to B\times B$ defined by $r_B\left(a,b\right)=\left(a\circ\left(\iota\left(a\right)+b\right),\overline{\left(\iota\left(a\right) +b\right)}\circ b \right)$ for all $a,b \in B$ is a solution if and only if the following condition holds
	\begin{align}\label{eq:condSolution}
	a + \lambda_b\left(c\right)\circ \left(\iota\left(1\right) +\rho_c\left(b\right)\right) = a + b \circ\left(\iota\left(1\right)+c\right)
	\end{align}
	for all $a,b,c \in B$.
	\begin{proof}
		It is easily verified that $\left(B,r_B\right)$ is a solution if and only if, for all $a,b,c \in B$,
		\begin{align*}
		\left(\lambda_a\lambda_b\left(c\right),\lambda_{\rho_{\lambda_b\left(c\right)}\left(a\right)}\rho_c\left(b\right),\rho_{\rho_c\left(b\right)}\rho_{\lambda_b\left(c\right)}\left(a\right)\right)=
		\left(\lambda_{\lambda_a\left(b\right)}\lambda_{\rho_b\left(a\right)}\left(c\right),\rho_{\lambda_{\rho_b\left(a\right)}\left(c\right)}\lambda_a\left(b\right),\rho_c\rho_b\left(a\right)\right).
		\end{align*}
		Denote the first triple by $\left(s_1, s_2, s_3\right)$ and the second by $\left(t_1, t_2, t_3\right)$. Since, for all $x,y \in B$,
		\begin{align}\label{eq:lambdarho}
		\lambda_x\left(y\right)\circ\rho_y\left(x\right)=x\circ y
		\end{align}
		holds, it follows that
		\begin{align*}
		s_1\circ s_2 \circ s_3 
		&= \lambda_a\lambda_b\left(c\right)\circ\lambda_{\rho_{\lambda_b\left(c\right)}\left(a\right)}\rho_c\left(b\right)\circ\rho_{\rho_c\left(b\right)}\rho_{\lambda_b\left(c\right)}\left(a\right) \\
		&= \lambda_a\lambda_b\left(c\right)\circ \rho_{\lambda_b\left(c\right)}\left(a\right)\circ \rho_c\left(b\right)\\
		&=a \circ \lambda_b\left(c\right)\circ \rho_c\left(b\right)\\
		&= a\circ b \circ c
		\end{align*}
		and
		\begin{align*}
		t_1\circ t_2 \circ t_3 &= \lambda_{\lambda_a\left(b\right)}\lambda_{\rho_b\left(a\right)}\left(c\right)\circ \rho_{\lambda_{\rho_b\left(a\right)}\left(c\right)}\lambda_a\left(b\right)\circ\rho_c\rho_b\left(a\right)\\
		&=\lambda_a\left(b\right) \circ \lambda_{\rho_b\left(a\right)}\left(c\right)\circ \rho_b\left(a\right)\\
		&=\lambda_a\left(b\right) \circ \rho_b\left(a\right) \circ c\\
		&= a \circ b \circ c.
		\end{align*}
		Thus, $s_1\circ s_2\circ s_3 = t_1 \circ t_2 \circ t_3$. \\
		Now, suppose that \eqref{eq:condSolution} holds. 
		As by \cref{lemma:lambda}, $\lambda:\left(B,\circ\right)\to \End\left(B,+\right)$ is a homomorphism, and since \eqref{eq:lambdarho} it follows that
		\begin{align}\label{eq:t1ands1}
		t_1 = \lambda_{\lambda_a\left(b\right)}\lambda_{\rho_b\left(a\right)}\left(c\right) = \lambda_{\lambda_a\left(b\right)\circ \rho_b\left(a\right)}\left(c\right) = \lambda_{a\circ b}\left(c\right) = \lambda_a\lambda_b\left(c\right)=s_1.
		\end{align}
		To prove that $s_3=t_3$, first note that for any $x,y,z\in B$, it holds
		\begin{align*}
		    \rhoo{z}{\rhoo{y}{\left(x\right)}}
		    &=\overline{\left(\iota\rhoo{y}{\left(x\right)}+z\right)}\circ z\\
		    &=\overline{\left(\iota\left(\overline{\left(\iota\left(x\right)+y\right)\circ y}\right) +z\right)}\circ z\\
		    &=\overline{\left(\bar{y}\circ\iota\left(\overline{\iota\left(x\right)+y}\right)+z\right)}\circ z &\mbox{(by \eqref{eq:defAlmost1})}\\
		    &=\overline{\left(\bar{y}\circ\left(\iota\left(x\right) +y\right)\circ\iota\left(1\right)+z \right)}\circ z &(\mbox{by \cref{lemma:24.2}})\\
		    &=\overline{\left(\bar{y}\circ\left(\iota\left(x\right)+y\circ\iota\left(1\right)\right) +z\right)}\circ z &\mbox{(by \eqref{eq:defAlmost2})}\\
		    &=\overline{\left(\bar{y}\circ\left(\iota\left(x\right)+y\circ\iota\left(1\right)\right) +\lambdaa{1}{\left(z\right)}\right)}\circ z&(\mbox{by \cref{lemma:24.1}})\\
		    &=\overline{\left(\bar{y}\circ\left(\iota\left(x\right)+y\circ\iota\left(1\right)\right) +\lambdaa{\bar{y}}{\lambdaa{y}{\left(z\right)}}\right)}\circ z&\mbox{(by \cref{lemma:lambda_almost})}\\
		    &=\overline{\bar{y}\circ\left(\iota\left(x\right)+y\circ\iota\left(1\right) +\lambdaa{y}{\left(z\right)}\right)}\circ z &\mbox{(by \eqref{eq:defAlmost3})}\\
		    &=\overline{\left(\iota\left(x\right)+y\circ\left(\iota\left(1\right)+z\right)\right)}\circ y\circ z &\mbox{(by \eqref{eq:defAlmost3})}
		\end{align*}
		Hence, by \eqref{eq:condSolution} and \eqref{eq:lambdarho}, it follows that
		\begin{align*}
		    s_3&=\rhoo{\rhoo{c}{\left(b\right)}}{\rhoo{\lambdaa{b}{\left(c\right)}}{\left(a\right)}}\\
		    &=\overline{\left(\iota\left(a\right)+\lambdaa{b}{\left(c\right)}\circ\left(\iota\left(1\right)+\rhoo{c}{\left(b\right)}\right)\right)}\circ \lambdaa{b}{\left(c\right)}\circ \rhoo{c}{\left(b\right)}\\
		    &=\overline{\left(\iota\left(a\right)+b\circ\left(\iota\left(1\right)+c\right)\right)}\circ b\circ c\\
		    &=t_3.
		\end{align*}
		Moreover, since $s_1\circ s_2\circ s_3=t_1\circ t_2\circ t_3$, $t_1=s_1$, $s_3=t_3$, and $\left(B,\circ\right)$ is a group it holds $s_2=t_2$.
		
		Conversely, suppose that $\left(B,r_B\right)$ is a solution. In particular, $$\overline{\left(\iota\left(a\right) +\lambda_b\left(c\right)\circ\left(\iota\left(1\right)+\rho_c\left(b\right)\right) \right)}\circ b\circ c=s_3=t_3 = \overline{\left(\iota\left(a\right) + b \circ \left(\iota\left(1\right)+c\right)\right)}\circ b \circ c.$$Hence \eqref{eq:condSolution} holds,  since, by \cref{lemma:24.2}, $\iota$ is bijective.
	\end{proof}
\end{theor}

Let $\left(B,+,\circ,\iota\right)$ be an almost left semi-brace that satisfies \eqref{eq:defAlmost2} and \eqref{eq:condSolution}. The map $r_B$ defined in \cref{theor:solutionNew} is said to be the \emph{solution associated to $B$}.

Note that if $\left(B,+,\circ, \iota\right)$ is an almost left semi-brace such that \eqref{eq:defAlmost2} holds and $\rho:B \to \Map\left(B,B\right)$ is an anti-homomorphism, then \eqref{eq:condSolution} is satisfied. Indeed, since $\iota$ is bijective, it holds
\begin{align*}
    a+ \lambdaa{b}{\left(c\right)}\circ\left(\iota\left(1\right)+\rhoo{c}{\left(b\right)}\right)
    &=\lambdaa{b}{\left(c\right)}\circ\rhoo{c}{\left(b\right)}\circ\overline{\rhoo{\rhoo{c}{\left(b\right)}}{\rhoo{\lambdaa{b}{\left(c\right)}}{\left(\iota^{-1}\left(a\right)\right)}}}\\
    &=b\circ c \circ\overline{\rhoo{\lambdaa{b}{\left(c\right)}\circ\rhoo{c}{\left(b\right)}}{\left(\iota^{-1}\left(a\right)\right)}}\\
    &=b\circ c \circ\overline{\rhoo{b\circ c}{\left(\iota^{-1}\left(a\right)\right)}}\\
    &=b\circ c \circ \overline{\rhoo{b}{\rhoo{c}{\left(\iota^{-1}\left(a\right)\right)}}}\\
    &= a + b \circ\left(\iota\left(1\right) +c\right),
\end{align*}
i.e., \eqref{eq:condSolution} is satisfied.

Following \cite[p.105]{CeJO14}, we define a homomorphism between set theoretical solutions.

\begin{defin}
	Let $X,Y$ be sets and let $r:X \times X\to X \times X$ and $r':Y\times Y \to Y \times Y$ be maps such that $r\left(a,b\right) = \left(\lambda_a\left(b\right),\rho_b\left(a\right)\right)$ and $r'\left(x,y\right) = \left(\lambda'_a\left(b\right),\rho'_b\left(a\right)\right)$. A \emph{homomorphism} from $\left(X,r\right)$ to $\left(Y,r'\right)$ is a map $f : X \to Y$ such that $\left(f\times f\right)r = \left(f \times f \right)r'$, i.e., $\left(f\left(\lambda_a\left(b\right)\right),f\left(\rho_b\left(a\right)\right)\right) = \left(\lambda'_{f\left(a\right)}\left(f\left(b\right)\right),\rho'_{f\left(b\right)}\left(f\left(a\right)\right)\right)$, for all $a,b \in X$. If $f$ is bijective then we say that $\left(X,r\right)$ and $\left(Y,r'\right)$ are \emph{isomorphic solutions}.
\end{defin}

In \cite[Theorem 3]{CCoSt19x}, it is proved that for any left semi-brace $B$, the map $r_B:B\times B \to B\times B$ defined by $r_B\left(a,b\right)=\left(a\circ\left(\overline{a}+b\right),\overline{\left(\overline{a}+b\right)}\circ b\right)$ is a solution if and only if 
\begin{align}\label{eq:condSolutionSemi}
a + \lambda_b\left(c\right)\circ \left(1+\rho_c\left(b\right)\right) = a + b \circ\left(1+c\right),
\end{align}
where $\lambda_b\left(c\right)=b\circ\left(\overline{b}+c\right)$, $\rho_c\left(b\right)=\overline{\left(\overline{b}+ c\right)}\circ c$, and $1$ is the identity of the group $\left(B,+\right)$.

In this theorem we prove that if $\left(B,+,\circ,\iota\right)$ satisfies \eqref{eq:condSolution}, then $\left(B,\oplus,\circ^{op}\right)$ satisfies \eqref{eq:condSolutionSemi} and that the solution associated to the almost left semi-brace and the solution associated to the left semi-brace are isomorphic.

\begin{theor}
	Let $\left(B,+,\circ, \iota\right)$ be an almost left semi-brace that satisfies \eqref{eq:defAlmost2} and \eqref{eq:condSolution}, and let $\left(B,\oplus,\circ^{op}\right)$ be the left semi-brace associated with $B$. Then $\left(B,\oplus,\circ^{op}\right)$ satisfies \eqref{eq:condSolutionSemi}. Furthermore, if
	$r_B$ is the solution defined in \cref{theor:solutionNew}, and $r'_B$ the solution associated to the left semi-brace $\left(B,\oplus,\circ^{op}\right)$ as defined in \cite[Theorem 3]{CCoSt19x}, then $r_B$ and $r'_B$ are isomorphic.
	\begin{proof}
		First, we compute the maps $\lambda'_a\left(b\right) =a\circ^{op}\left(\overline{a}\oplus b\right)$
		and $\rho'_b\left(a\right)=\overline{\left(\overline{a}\oplus b\right)}\circ^{op}b$, for the left semi-brace $\left(B,\oplus,\circ^{op}\right)$. Let $a,b,c\in B$. By \cref{lemma:24.2}, \eqref{eq:iotalemma}, and \eqref{eq:defAlmost2}, we have 
		\begin{align*}
		\lambda'_a\left(b\right) 
		&=a \circ^{op}\left(\bar{a}\oplus b\right)
		= \iota^{-1}\left(\iota\left(\bar{a}\right)+\iota\left(b\right)\right)\circ a \\
		&= \iota^{-1}\left(\bar{a}\circ \left(\iota\left(\bar{a}\right)+\iota\left(b\right)\right)\right)\\
		&=\iota^{-1}\left(\bar{a}\circ\left(\iota\left(\bar{a}\right)+b\circ\iota\left(1\right)\right)\right)\\
		&= \iota^{-1}\left(\bar{a}\circ \left(\iota\left(\bar{a}\right) + \bar{b}\right)\circ \iota\left(1\right)\right)\\
		&=1 \circ \overline{\left(\bar{a}\circ \left(\iota\left(\bar{a}\right) + \bar{b}\right)\right)}\\
		&=\overline{\lambda_{\bar{a}}\left(\bar{b}\right)}.
		\end{align*}
		Moreover, by \cref{lemma:24.2} we obtain $\bar{b} = \iota\left(b\right)\circ\overline{\iota\left(1\right)}$. Using also \eqref{eq:defAlmost2}, it follows that
		\begin{align*}
		\rho'_b\left(a\right) &= \overline{\left(\bar{a}\oplus b\right)}\circ^{op} b 
		= b \circ \overline{\iota^{-1}\left(\iota\left(\bar{a}\right)+\iota\left(b\right)\right)}\\
		&= b \circ \iota\iota^{-1}\left(\iota\left(\bar{a}\right) +\iota\left(b\right)\right)\circ \overline{\iota\left(1\right)}\\
		&=b \circ \left(\iota\left(\bar{a}\right)+\bar{b}\circ\iota\left(1\right)\right)\circ\overline{\iota\left(1\right)}\\
		&=b \circ\left(\iota\left(\bar{a}\right)+\bar{b}\right)\circ\iota\left(1\right)\circ\overline{\iota\left(1\right)}\\
		&=b \circ \left(\iota\left(\bar{a}\right) +\bar{b}\right)\\
		&= \overline{\overline{\left(\iota\left(\bar{a}\right) +\bar{b}\right)}\circ \bar{b}}\\
		&= \overline{\rho_{\bar{b}}\left(\bar{a}\right)}.
		\end{align*}
		Now, suppose that $\left(B,+,\circ,\iota\right)$ satisfies \eqref{eq:defAlmost2} and \eqref{eq:condSolution}.
		By \cref{lemma:24.2} and by \eqref{eq:defAlmost2}, we have that
		\begin{align*}
		    \iota\left(\overline{a+b}\right) = \left(\iota\iota^{-1}\left(a\right)+b\right)\circ\iota\left(1\right)
		    =\iota\iota^{-1}\left(a\right)+b\circ\iota\left(1\right) = a +\iota\left(\bar{b}\right).
		\end{align*}
		And using also \eqref{eq:defAlmost1}, it follows that
		\begin{align*}
		a\oplus b \circ^{op}\left(1\oplus c\right)
		&=\iota^{-1}\left(\iota\left(a\right)+\iota\left(\iota^{-1}\left(\iota\left(1\right)+\iota\left(c\right)\right)\circ b\right)\right)\\
		&=\iota^{-1}\left(\iota\left(a\right)+\overline{b}\circ\left(\iota\left(1\right)+\iota\left(c\right)\right)\right)\\
		&=\iota^{-1}\left(\iota\left(a\right)+\overline{b}\circ\left(\iota\left(1\right)+\overline{c}\circ \iota\left(1\right)\right)\right)\\
		&=\iota^{-1}\left(\iota\left(a\right)+\overline{b}\circ\left(\iota\left(1\right)+\overline{c}\right)\circ\iota\left(1\right)\right)\\
		&=\iota^{-1}\left(\left(\iota\left(a\right)+\overline{b}\circ\left(\iota\left(1\right)+\overline{c}\right)\right)\circ\iota\left(1\right)\right)
		\end{align*}
		and, by \eqref{eq:condSolution},
		\begin{align*}
		a\oplus\lambda'_a\left(b\right)\circ^{op}\left(1\oplus\rho'_c\left(b\right)\right)
		&=\iota^{-1}\left(\iota\left(a\right) +\iota\left(\iota^{-1}\left(\iota\left(1\right) + \iota\left(\overline{\rho_{\overline{c}}\left(\overline{b}\right)}\right) \right)\circ\overline{\lambda_{\overline{b}}\left(\overline{c}\right)} \right) \right)\\
		&=\iota^{-1}\left( \iota\left(a\right) + \lambda_{\overline{b}}\left(\overline{c}\right)\circ\left(\iota\left(1\right) + \rho_{\overline{c}}\left(\overline{b}\right)\circ \iota\left(1\right) \right) \right) \\
		&=\iota^{-1}\left( \iota\left(a\right) + \lambda_{\overline{b}}\left(\overline{c}\right)\circ\left(\iota\left(1\right) + \rho_{\overline{c}}\left(\overline{b}\right) \right) \circ \iota\left(1\right)\right) \\
		&=\iota^{-1}\left( \left(\iota\left(a\right) + \lambda_{\overline{b}}\left(\overline{c}\right)\circ\left(\iota\left(1\right) + \rho_{\overline{c}}\left(\overline{b}\right) \right)  \right) \circ\iota\left(1\right)\right)\\
		&=\iota^{-1}\left(\left(\left(\iota\left(a\right)+\overline{b}\circ\left(\iota\left(1\right)+\overline{c}\right)\right)\right)\circ\iota\left(1\right)\right)\\
		&=\iota^{-1}\left(\left(\iota\left(a\right)+\overline{b}\circ\left(\iota\left(1\right)+\overline{c}\right)\right)\circ\iota\left(1\right)\right),
		\end{align*}
		i.e., \eqref{eq:condSolutionSemi} holds in the left semi-brace $\left(B,\oplus,\circ^{op}\right)$.
		Therefore, the map $r'_B\left(a,b\right) := \left(\lambda'_a\left(b\right), \rho'_b\left(a\right)\right)$ is a solution.
		Finally, defining $f:B\to B$ by $f\left(a\right)=\bar{a}$, it holds that
		\begin{align*}
		r'_B\left(f\times f\right)\left(a,b\right)
		&=r'_B\left(\bar{a},\bar{b}\right)
		=\left(\lambda'_{\bar{a}}\left(\bar{b}\right),\rho'_{\bar{b}}\left(\bar{a}\right)\right)\\
		&=\left(\overline{\lambda_a\left(b\right)},\overline{\rho_b\left(a\right)}\right)= \left(f\times f \right)\left(\lambda_a\left(b\right),\rho_b\left(a\right)\right)\\
		&=\left(f\times f\right) r_B\left(a,b\right),
		\end{align*}
		for all $a,b \in B$, i.e., $r'_B$ and $r_B$ are isomorphic solutions.
	\end{proof}
\end{theor}

\subsection*{Acknowledgments}
The first author is member of GNSAGA (INdAM). The second author is supported by Fonds voor Wetenschappelijk Onderzoek (Flanders), via an FWO Aspirant-mandate.

\bibliographystyle{abbrv}
	\bibliography{bibliography}

\def\cprime{$'$}
\begin{thebibliography}{10}

\bibitem{BaCeJeOk19}
D.~Bachiller, F.~Ced\'{o}, E.~Jespers, and J.~Okni\'{n}ski.
\newblock Asymmetric product of left braces and simplicity; new solutions of
  the {Y}ang--{B}axter equation.
\newblock {\em Commun. Contemp. Math.}, 21(8):1850042, 30, 2019.

\bibitem{Br18}
T.~Brzezi\'{n}ski.
\newblock Towards semi-trusses.
\newblock {\em Rev. Roumaine Math. Pures Appl.}, 63(2):75--89, 2018.

\bibitem{Br19}
T.~Brzezi\'{n}ski.
\newblock Trusses: between braces and rings.
\newblock {\em Trans. Amer. Math. Soc.}, 372(6):4149--4176, 2019.

\bibitem{Br20}
T.~Brzezi{\'n}ski.
\newblock Trusses: Paragons, ideals and modules.
\newblock {\em Journal of Pure and Applied Algebra}, 224(6):106258, 2020.

\bibitem{Br19xa}
T.~Brzezi{\'n}ski and B.~Rybo{\l}owicz.
\newblock Modules over trusses vs modules over rings: direct sums and free
  modules.
\newblock {\em arXiv:1909.05807}, 2019.

\bibitem{Br19xb}
T.~Brzezi{\'n}ski and B.~Rybo{\l}owicz.
\newblock On congruence classes and extensions of rings with applications to
  braces.
\newblock {\em arXiv:1912.00907}, 2019.

\bibitem{CCoSt17}
F.~Catino, I.~Colazzo, and P.~Stefanelli.
\newblock Semi-braces and the {Y}ang-{B}axter equation.
\newblock {\em J. Algebra}, 483:163--187, 2017.

\bibitem{CCoSt20}
F.~Catino, I.~Colazzo, and P.~Stefanelli.
\newblock Algebraic tools for {S}olving the {Y}ang--{B}axter equation:
  {G}eneralized semi-braces and the {Y}ang--{B}axter equation.
\newblock {\em Oberwolfach Reports}, 2019(51), 2019.

\bibitem{CCoSt19x}
F.~Catino, I.~Colazzo, and P.~Stefanelli.
\newblock The matched product of the solutions to the {Y}ang--{B}axter equation
  of finite order.
\newblock {\em arXiv preprint arXiv:1904.07557}, 2019.

\bibitem{CCoSt19}
F.~Catino, I.~Colazzo, and P.~Stefanelli.
\newblock Skew left braces with non-trivial annihilator.
\newblock {\em J. Algebra Appl.}, 18(2):1950033, 23, 2019.

\bibitem{CeGaIvSm18}
F.~Ced\'{o}, T.~Gateva-Ivanova, and A.~Smoktunowicz.
\newblock Braces and symmetric groups with special conditions.
\newblock {\em J. Pure Appl. Algebra}, 222(12):3877--3890, 2018.

\bibitem{CeJO14}
F.~Ced\'{o}, E.~Jespers, and J.~Okni\'{n}ski.
\newblock Braces and the {Y}ang--{B}axter equation.
\newblock {\em Comm. Math. Phys.}, 327(1):101--116, 2014.

\bibitem{CeJeVe20}
F.~Ced\'o, E.~Jespers, and C.~Verwimp.
\newblock Structure monoids of set-theoretic solutions of the {Y}ang--{B}axter
  equation.
\newblock {\em arXiv:1912.09710}, 2019.

\bibitem{CSV}
F.~Ced\'{o}, A.~Smoktunowicz, and L.~Vendramin.
\newblock Skew left braces of nilpotent type.
\newblock {\em Proc. Lond. Math. Soc. (3)}, 118(6):1367--1392, 2019.

\bibitem{Drinfeld}
V.~G. Drinfel{\cprime}d.
\newblock On some unsolved problems in quantum group theory.
\newblock In {\em Quantum groups ({L}eningrad, 1990)}, volume 1510 of {\em
  Lecture Notes in Math.}, pages 1--8. Springer, Berlin, 1992.

\bibitem{GV17}
L.~Guarnieri and L.~Vendramin.
\newblock Skew braces and the {Y}ang--{B}axter equation.
\newblock {\em Math. Comp.}, 86(307):2519--2534, 2017.

\bibitem{JeKuVAV19}
E.~Jespers, L.~Kubat, A.~Van~Antwerpen, and L.~Vendramin.
\newblock Factorizations of skew braces.
\newblock {\em Math. Ann.}, 375(3-4):1649--1663, 2019.

\bibitem{JeKuVAV20}
E.~Jespers, {\L}.~Kubat, A.~Van~Antwerpen, and L.~Vendramin.
\newblock Radical and weight of skew braces and their applications to structure
  groups of solutions of the {Y}ang--{B}axter equation.
\newblock {\em arXiv:2001.10967}, 2020.

\bibitem{JVA19}
E.~Jespers and A.~Van~Antwerpen.
\newblock Left semi-braces and solutions of the {Y}ang--{B}axter equation.
\newblock {\em Forum Math.}, 31(1):241--263, 2019.

\bibitem{konovalov2018skew}
A.~Konovalov, A.~Smoktunowicz, and L.~Vendramin.
\newblock On skew braces and their ideals.
\newblock {\em Experimental Mathematics}, pages 1--10, 2018.

\bibitem{LeVe16}
V.~Lebed and L.~Vendramin.
\newblock Cohomology and extensions of braces.
\newblock {\em Pacific J. Math.}, 284(1):191--212, 2016.

\bibitem{Mi18}
M.~M. Miccoli.
\newblock Almost semi-braces and the {Y}ang--{B}axter equation.
\newblock {\em Note Mat.}, 38(1):83--88, 2018.

\bibitem{Ru07}
W.~Rump.
\newblock Braces, radical rings, and the quantum {Y}ang--{B}axter equation.
\newblock {\em J. Algebra}, 307(1):153--170, 2007.

\end{thebibliography}

\end{document}